\documentclass[12pt]{amsart}
\usepackage{amsfonts,amsmath,amsthm,amssymb,url}
\usepackage[russian]{babel}
\usepackage[utf8]{inputenc}

\renewcommand{\le}{\leqslant}
\renewcommand{\ge}{\geqslant}
\renewcommand{\leq}{\leqslant}
\renewcommand{\geq}{\geqslant}
\def\R{{\mathbb R}} \def\Z{{\mathbb Z}} \def\F{{\mathbb F}}

\newtheorem{theorem}{Теорема}
    \newtheorem{lemma}[theorem]{Лемма}
    
    \newtheorem{proposition}[theorem]{Утверждение}

\usepackage{hyperref}
\hypersetup{
   colorlinks   = true, 
   urlcolor     = blue, 
   linkcolor    = blue, 
   citecolor   = red 
}

\newtheoremstyle{mydefinition}
  {3pt}
  {3pt}
  {\normalfont}
  {\parindent}
  {\bfseries}
  {.}
  { }
  {}

\theoremstyle{mydefinition}
\newtheorem{pr}[theorem]{Задача}

\begin{document}

\newcommand{\mytitle}{Мотивированное изложение комбинаторной теоремы о нулях}


\title{\mytitle}
\author{М. Ложкин и А. Скопенков}
\thanks{М. Ложкин: \texttt{lozhkin.mixail@gmail.com}.
Национальный Исследовательский Университет <<Высшая Школа Экономики>>.
\newline
А. Скопенков: \texttt{https://users.mccme.ru/skopenko}; \texttt{skopenko@mccme.ru}.
Московский Физико-Технический Институт, Независимый Московский Университет.
\newline
Заметка возникла из обсуждений семинаров по курсу <<Дискретный анализ>> А. Райгородского на ФПМИ МФТИ, а также занятий кружка ЦПМ, проводимых М. Ложкиным.
Благодарим А. Волостнова за предосталение его материалов по \S2, а также А. Волостнова, Д. Гринберга, Н. Ленскую и Ф. Петрова и анонимного рецензента за полезные обсуждения.}

\date{}
\maketitle

О применениях линейной алгебры и многочленов в комбинаторике хорошо известно (см., например, [R1, R2, S]).
Менее известные в школьном образовании применения \emph{комбинаторной теоремы о нулях} (теорем Алона \ref{c:degb}.b и \ref{c:degb'}.b) описаны в [A, C, D, KS].
Мы покажем, как эта теорема
постепенно и естественно возникает при решении <<олимпиадных>> задач.
Простейшие идеи (утверждения \ref{l:al2}, \ref{l:al2'}, \ref{l:alp}, \ref{l:alp'}, \ref{c:deg0p} и \ref{c:deg0p'}) естественно появляются при решении задач
\ref{p:poly2}, \ref{p:poly} и \ref{p:circ}, чем подводят читателя к обобщениям (теоремам \ref{c:degb} и \ref{c:degb'}).
Явная формулировка промежуточных частных случаев позволяет читателю проделать самостоятельно возможно большее количество этих постепенных обобщений.


Многие излагаемые здесь результаты близки к переднему краю науки.
Однако эти результаты <<олимпиадные>>: для понимания формулировок и придумывания доказательств не требуется знаний, выходящих за пределы <<кружковской>> программы (или программы первого курса).

Подумайте самостоятельно над доказательством каждого утверждения!
Доказательства приведены после формулировки или в конце заметки.
Если условие задачи является утверждением, то задача состоит в том, чтобы это утверждение доказать
(тогда мы пишем <<доказательство утверждения>>, а не <<решение задачи>>).

\section{Приходим к комбинаторной теореме о нулях}

\emph{Степенью} одночлена $x_1^{d_1}\ldots x_n^{d_n}$ называется число $d_1+\ldots+d_n$.
\emph{Степенью} $\deg f$ ненулевого многочлена $f$ называется максимальная степень входящих в него одночленов.
\emph{Степенью} нулевого многочлена называется минус бесконечность.

\begin{pr}\label{p:poly2} Даны два многочлена по модулю 2 от $n$ переменных, сумма степеней которых меньше $n$.
Если они имеют общий корень, то они имеют хотя бы два общих корня.
\end{pr}

Доказательство основано на следующем утверждении \ref{l:al2}.

Через $\Z_p$ обозначается множество вычетов по модулю $p$.
Напомним, что
$$A_1\times\ldots\times A_n\ :=\ \{(a_1,\ldots,a_n)\ :\ a_1\in A_1,\ldots,a_n\in A_n\}.$$
В утверждениях
\ref{l:al2} и \ref{l:al2'}
сумирование и равенства рассматриваются по модулю 2.

\begin{proposition}\label{l:al2} Дан многочлен $f$ по модулю 2 от $n$ переменных степени меньше $n$.

(a) Имеем $\sum\limits_{\alpha\in\Z_2^n} f(\alpha) = 0.$

(b) Если $f$ имеет корень, то он имеет хотя бы два корня.
\end{proposition}

Назовем многочлен (с какими угодно коэффициентами) от нескольких переменных {\bf неполным}, если в каждом его одночлене одна из переменных отсутствует.

\smallskip
{\it Доказательство утверждения \ref{l:al2}.}
(a) Так как $\deg f<n$, то $f$ неполный.
Если $f$ не содержит $x_1$, то  $f(0,\beta)=f(1,\beta)$ при любом $\beta\in\Z_2^{n-1}$, поэтому равенство из утверждения верно.
Значит, оно верно для любого неполного многочлена $f$.

(b) Применим п. (a).
Так как сумма четного количества слагаемых равна нулю и одно из них равно нулю, то еще одно равно нулю.
\qed

\smallskip
{\it Доказательство утверждения \ref{p:poly2}.}
Обозначим данные многочлены через $f_1,f_2$.
Обозначим $g := f_1f_2+f_1+f_2 = 1+(1+f_1)(1+f_2)$.
Набор $\alpha\in\Z_2^n$ является корнем многочлена $g$ тогда и только тогда, когда $\alpha$ является общим корнем многочленов $f_1,f_2$.
Имеем $\deg g \le \deg f_1 + \deg f_2 < n$.
Значит, по утверждению \ref{l:al2}.b многочлены $f_1,f_2$ имеют хотя бы два общих корня.
\qed

{\renewcommand{\thetheorem}{\ref{l:al2}$'$}\addtocounter{theorem}{-1}
\begin{proposition}\label{l:al2'} Дан многочлен $f$ по модулю 2 от $n$ переменных степени не более $n$.

(a) Коэффициент при одночлене $x_1\cdots x_n$ равен
$\sum\limits_{\alpha\in\Z_2^n} f(\alpha).$

(b) Если $f$ зануляется в каждой точке, то его коэффициент при одночлене $x_1\cdots x_n$ равен нулю.
\end{proposition}
}

{\it Доказательство.}
(a) Утверждение \ref{l:al2}.a доказано для неполных многочленов (не обязательно имеющих степень менее $n$).
Поэтому достаточно доказать утверждение \ref{l:al2'}.a для $f(x_1,\ldots,x_n) = kx_1\cdots x_n$.
Для него сумма из условия равна $k$.


(b) Следует из п. (a).
\qed

\smallskip
Далее число $p$ простое.
В утверждениях \ref{p:poly}, \ref{l:alp}, \ref{l:alp'}, \ref{c:deg0p} и \ref{c:deg0p'}
сумирование и равенства рассматриваются по модулю $p$.

\begin{pr}[теорема Шевалле]\label{p:poly}
Даны многочлены по простому модулю $p$ от $n$ переменных, сумма степеней которых меньше $n$.
Если они имеют общий корень, то они имеют хотя бы два общих корня.
\end{pr}

Доказательство основано на следующем естественном обобщении утверждения \ref{l:al2}.
Из него вытекает более сильный факт: количество корней делится на $p$ (теорема Варнинга).

\begin{proposition}\label{l:alp} Дан многочлен $f$ по модулю $p$ от $n$ переменных степени меньше
\linebreak
$n(p-1)$.

(a) Имеем $\sum\limits_{\alpha\in\Z_p^n} f(\alpha) = 0.$

(b) Если $f$ не зануляется хотя бы в одной точке, то он не зануляется хотя бы в двух точках.
\end{proposition}

{\it Доказательство.}
(a) Ввиду линейности достаточно доказать утверждение для $f(x_1,\ldots,x_n) = x_1^{k_1}\cdots x_n^{k_n}$.
В этом случае выражение в формуле можно разложить на множители:
$$
\sum\limits_{\alpha\in\Z_p^n} f(\alpha) = \sum_{\alpha_1\in\Z_p}\alpha_1^{k_1}\cdot\ldots\cdot\sum_{\alpha_n\in\Z_p}\alpha_n^{k_n}.
$$
Так как $k_1+\ldots+k_n<n(p-1)$, то, не уменьшая общности, $k_1<p-1$.
Тогда\footnote{Равенство $\sum_{\alpha_1\in\Z_p}\alpha_1^{k_1}=0$ следует также из леммы \ref{p:lagrk}.b для $A=\Z_p$ и $m=k_1$.}, беря первообразный корень $g$ по модулю $p$, получаем:
$$
\sum_{\alpha_1\in\Z_p} \alpha_1^{k_1} = \sum_{0\le n<p-1} g^{nk_1} = \frac{g^{k_1(p-1)}-1}{g^{k_1}-1} = 0,\quad\text{где}
$$

$\bullet$ первое равенство верно, поскольку $g$~---~первообразный корень;

$\bullet$ второе равенство верно, поскольку $g^{k_1}\ne1$;

$\bullet$ последнее равенство верно, поскольку по малой теореме Ферма $g^{k_1(p-1)}=1$.

(b) Применим п. (a).
Так как сумма равна нулю и одно из слагаемых не равно нулю, то еще одно не равно нулю.
\qed


\smallskip
{\it Доказательство утверждения \ref{p:poly}.}
Обозначим данные многочлены через $f_1,f_2,\ldots,f_k$.
Обозначим
\[
g := \prod_{j = 1}^k (f_j^{p - 1} - 1).
\]
Ввиду малой теоремы Ферма, $g(\alpha)\ne0$ тогда и только тогда, когда $\alpha$ является общим корнем многочленов $f_1,\ldots,f_k$.
Так как многочлены $f_1,\ldots,f_k$ имеют общий корень, то многочлен $g$ не зануляется в некоторой точке.
Имеем
\[
\deg g = (p - 1)\left(\deg f_1 + \deg f_2 + \ldots + \deg f_k\right) < n(p - 1).
\]
Значит, по утверждению \ref{l:alp}.b многочлен $g$ не зануляется хотя бы в двух точках.
Поэтому у многочленов $f_1,\ldots,f_k$ есть хотя бы два общих корня.
\qed

{\renewcommand{\thetheorem}{\ref{l:alp}$'$}\addtocounter{theorem}{-1}
\begin{proposition}\label{l:alp'} Дан многочлен $f$ по модулю $p$ от $n$ переменных степени не более $n(p-1)$.

(a) Коэффициент при одночлене $x_1^{p-1}\cdots x_n^{p-1}$ равен $(-1)^n\sum\limits_{\alpha\in\Z_p^n} f(\alpha)$.

(b) Если $f$ зануляется в каждой точке, то его коэффициент при одночлене $x_1^{p-1}\cdots x_n^{p-1}$ равен нулю.
\end{proposition}
}

{\it Доказательство.}
(a) Утверждение \ref{l:alp}.a доказано для одночленов, у которых степень вхождения хотя бы одной переменной меньше $p-1$ (но не обязательно имеющих степень менее $n(p-1)$).
Поэтому достаточно доказать утверждение \ref{l:alp'}.a для
$f(x_1,\ldots,x_n) = kx_1^{p-1}\cdots x_n^{p-1}$.
Для него по малой теореме Ферма сумма из условия равна $(-1)^n(p-1)^nk=k$.

(b) Следует из п. (a).
\qed

\begin{pr}[Всероссийская олимпиада 2007/11.5]\label{p:circ}
В каждой вершине
100-угольника записаны два различных числа.
Тогда из каждой вершины можно удалить одно число так, чтобы оставшиеся числа в любых двух соседних вершинах  были различны.
\end{pr}

Доказательство основано на следующем естественном обобщении вышеприведенных утверждений \ref{l:alp'}, в котором вместо значений многочлена на $\Z_2^n$ и $\Z_p^n$ рассматриваются значения многочлена на произведении $A_1\times\ldots\times A_n\subset\Z_p^n$.

\smallskip
\emph{Замечание.}
У утверждений \ref{p:circ} и \ref{t:cada} имеется несложное прямое комбинаторное доказательство (см. доказательство утверждения \ref{t:cada} в решениях).
Однако для большей части утверждений этого текста доказательство, не использующее комбинаторную теорему о нулях, неизвестно (или сложно).

\smallskip
Далее в этом тексте $\F$ есть $\Z_p$ или $\R$ (случай произвольного \emph{поля}
доказывается так же, но не применяется в комбинаторных задачах из \S2).


\begin{proposition}\label{c:deg0p} Даны двухэлементные подмножества $A_1,\ldots,A_n\subset\F$ и многочлен $f$ с коэффициентами в $\F$ от $n$ переменных степени меньше $n$.

(a) Обозначим $\{a_{i0},a_{i1}\}:=A_i$.
Тогда
$$\sum\limits_{\alpha\in\Z_2^n}(-1)^{\alpha(1)+\ldots+\alpha(n)}f(a_{1\alpha(1)},\ldots,a_{n\alpha(n)}) = 0.$$


(b) Если $f$ не зануляется хотя бы в одной точке множества $A_1\times\ldots\times A_n$, то он не зануляется хотя бы в двух его точках.
\end{proposition}


{\it Доказательство.}
(a) Так как $\deg f<n$, то $f$ неполный.
Если $f$ не содержит $x_1$, то $f(a_{10},a_{2\alpha(2)},\ldots,a_{n\alpha(n)})=f(a_{11},a_{2\alpha(2)},\ldots,a_{n\alpha(n)})$
для любых $\alpha(2),\ldots,\alpha(n)\in\Z_2$, поэтому нужное равенство верно.
Тогда оно верно и для любого неполного многочлена $f$.

(b) Следует из п. (a).
\qed

\smallskip
С помощью утверждения \ref{c:deg0p} доказываются утверждения \ref{p:plan}.b и \ref{p:reg}.


{\renewcommand{\thetheorem}{\ref{c:deg0p}$'$}\addtocounter{theorem}{-1}
\begin{proposition}\label{c:deg0p'} Даны двухэлементные подмножества $A_1,\ldots,A_n\subset\F$ и многочлен $f$ с
коэффициентами в $\F$ от $n$ переменных степени не более $n$.

(a) Обозначим $\{a_{i0},a_{i1}\}:=A_i$.
Тогда
коэффициент при одночлене $x_1\cdots x_n$ равен
$$\frac{\sum\limits_{\alpha\in\Z_2^n} (-1)^{\alpha(1)+\ldots+\alpha(n)} f(a_{1\alpha(1)},\ldots,a_{n\alpha(n)})}
{\prod\limits_{i=1}^n(a_{i,0}-a_{i,1})}.$$

(b) Если $f$ зануляется на $A_1\times\ldots\times A_n$, то коэффициент при одночлене $x_1\cdots x_n$ равен нулю.
\end{proposition}
}

{\it Доказательство.} (a) Утверждение \ref{c:deg0p}.a доказано для неполных многочленов (не обязательно имеющих степень менее $n$).
Поэтому достаточно доказать утверждение  \ref{c:deg0p'}.a для $f(x_1,\ldots,x_n) = kx_1\cdots x_n$.
Имеем
$$\sum\limits_{\alpha\in\Z_2^n} (-1)^{\alpha(1)+\ldots+\alpha(n)} a_{1\alpha(1)}\cdot\ldots\cdot a_{n\alpha(n)} =
\prod\limits_{i=1}^n(a_{i,0}-a_{i,1}).$$
Поэтому сумма из условия равна $k$.

(b) Следует из п. (a).
\qed

\smallskip
{\it Доказательство утверждения \ref{p:circ}.}
Обозначим $n=100$.
Определим многочлен
\[
f(x_1,\ldots,x_n) = \prod_{i=1}^n(x_i-x_{i+1}),\quad\text{где}\quad x_{n+1} = x_1.
\]
Имеем $\deg f=n$.
Так как $n$ четно, то в многочлене $f$ коэффициент при одночлене $x_1\cdots x_n$ равен $2\ne0$.
Применим утверждение \ref{c:deg0p'}.b для множеств чисел, записанных в вершинах. 
Получим, что найдутся такие числа $a_1,\ldots,a_n$, записанные в вершинах, что $f(a_1,\ldots,a_n)\ne 0$.
Это и есть искомые числа.
\qed

\begin{theorem}[Алон]\label{c:degb} Даны конечные подмножества $A_1,\ldots,A_n \subset \F$
и многочлен $f$ с коэффициентами в $\F$ от $n$ переменных степени менее $|A_1|+\ldots+|A_n|-n$.

(a) Тогда существуют (не зависящие от $f$) $P(\alpha)\in\F\setminus\{0\}$, где $\alpha\in A_1\times\ldots\times A_n$, для которых
$$
\sum_{\alpha\in A_1\times\ldots\times A_n}\frac{f(\alpha)}{P(\alpha)}=0.
$$


(b) Если $f$ не зануляется хотя бы в одной точке множества $A_1\times\ldots\times A_n$, то он не зануляется
хотя бы в двух его точках.
\end{theorem}

П. (b) следует из п. (a).
Для доказательства теорем \ref{c:degb}.a и \ref{c:degb'}.a требуется лемма \ref{p:lagrk}.b.

{\renewcommand{\thetheorem}{\ref{c:degb}$'$}\addtocounter{theorem}{-1}
\begin{theorem}[Алон]\label{c:degb'} Даны конечные непустые подмножества $A_1,\ldots,A_n\subset\F$ и
многочлен $f$ с коэффициентами в $\F$ от $n$ переменных степени не более $|A_1|+\ldots+|A_n|-n$.

(a) Тогда существуют (не зависящие от $f$) $P(\alpha)\in\F\setminus\{0\}$, где $\alpha\in A_1\times\ldots\times A_n$,
для которых коэффициент при одночлене $x_1^{|A_1|-1}\cdots x_n^{|A_n|-1}$ равен
$$
\sum_{\alpha\in A_1\times\ldots\times A_n}\frac{f(\alpha)}{P(\alpha)}.
$$


(b) Если $f$ зануляется
на $A_1\times\ldots\times A_n$,
то коэффициент при одночлене $x_1^{|A_1|-1}\cdots x_n^{|A_n|-1}$
равен нулю.
\end{theorem}
}


Далее $A\subset F$ --- любое конечное подмножество.
Для $a\in A$ обозначим
\[
Aa = \prod\limits_{b\in A,\ b\ne a}(a-b).
\]

\begin{lemma}\label{p:lagrk} (a) (Интерполяционная формула Лагранжа)
Для любых подмножества $A\subset\F$ и многочлена $f$ с коэффициентами в $\F$ от одной переменной степени менее $|A|$ имеем
\[f(x) = \sum_{a\in A}\frac{f(a)\prod\limits_{b\in A,\ b\ne a}(x-b)}{Aa}.\]

(b) Для любого подмножества $A\subset\F$ имеем
\[
\sum_{a\in A}\frac{a^m}{Aa} = \begin{cases} 0,& m<|A|-1;\\ 1,& m = |A|-1. \end{cases}
\]
\end{lemma}

{\it Доказательство.}
(a) В обеих частях равенства написаны многочлены степени менее $|A|$.
Они совпадают во всех точках множества $A$, следовательно они равны.
(Заметим, что $|A|<p+1$ при $\F=\Z_p$.)

(b) Подставим в п. (a) многочлен $f(x) = x^m$ при $m<|A|$:
\[x^m = \sum_{a\in A}\frac{a^m\prod_{b\in A\setminus\{a\}}(x-b)}{Aa}.\]
Приравняв коэффициенты обеих частей при $x^{|A|-1}$, получаем требуемое равенство.
\qed

\smallskip
{\it Доказательство теорем \ref{c:degb}.a и \ref{c:degb'}.a.}
Положим $P(\alpha):=A_1\alpha_1\cdot\ldots\cdot A_n\alpha_n$.
Ввиду линейности достаточно доказать каждую из теорем для $f(x_1, \ldots, x_n) = x_1^{k_1}\cdots x_n^{k_n}$.
В этом случае выражение в формуле можно разложить на множители:
$$
\sum_{\alpha_1\in A_1} \cdots \sum_{\alpha_n\in A_n} \frac{\alpha_1^{k_1}\cdots\alpha_n^{k_n}}{P(\alpha)} =
\sum_{\alpha_1\in A_1} \frac{\alpha_1^{k_1}}{A_1\alpha_1} \cdot \ldots \cdot \sum_{\alpha_n\in A_n} \frac{\alpha_n^{k_n}}{A_n\alpha_n}.
$$
Используем формулу \ref{p:lagrk}.b.
Если $k_i<|A_i|-1$ для некоторого $i$, то один из множителей равен нулю.
Иначе $k_i=|A_i|-1$ для любого $i$, тогда каждый из множителей равен~$1$.
\qed

\smallskip
\emph{Замечания.}
(a) Приведённое доказательство теорем Алона \ref{c:degb} и \ref{c:degb'} отлично от имеющегося в [A, C, D].

(b) Теоремы Алона являются естественным обобщением интерполяционной формулы Лагранжа \ref{p:lagrk}.a.
В [KP] это подмечено, и теорема Алона \ref{c:degb'}.a применена к вычислению коэффициентов некоторого многочлена.
В [P] обобщается следующее <<качественное следствие>> интерполяционной формулы Лагранжа \ref{p:lagrk}.a: \emph{многочлен степени $n$ от одной переменной задается своими значениями в $n+1$ точке} (ср. теорему \ref{c:degb'}.b); обобщение применяется к комбинаторным \emph{тождествам}.

(c) Многочлен $f$ от переменных $x_1,\ldots,x_n$ назовем \emph{$(d_1,\ldots,d_n)$-многочленом}, если среди его одночленов $x_1^{k_1}\cdots x_n^{k_n}$, отличных от $x_1^{d_1}\cdots x_n^{d_n}$, нет таких, что $k_i\ge d_i$ для всех $i\in\{1,\ldots,n\}$.
В теоремах Алона \ref{c:degb} и \ref{c:degb'} обозначим $d_i = |A_i|-1$ для $i=1,\ldots,n$.
Тогда условия на степень многочлена можно ослабить соответственно до следующих:

$\bullet$ $f$ является $(d_1,\ldots,d_n)$-многочленом и коэффициент при $x_1^{d_1}\cdots x_n^{d_n}$ равен нулю;

$\bullet$ $f$ является $(d_1,\ldots,d_n)$-многочленом.

(d) Утверждение \ref{p:lagrk}.b можно доказать напрямую, используя разложение по строке определителя Вандермонда (см. определение в решении задачи \ref{p:snevily}.b).


\section{Применяем комбинаторную теорему о нулях}

Здесь мы приводим  комбинаторные применения.
Об алгебраических обобщениях теорем Алона, Шевалле и Варнинга (теоремы \ref{c:degb}, \ref{c:degb'} и утверждение \ref{p:poly}) в разных направлениях см., например, [BS, C14, CGS].

\begin{pr}[теорема Эрдеша-Гинзбурга-Зива]\label{p:setp}
Для всякого простого $p$ любой набор из $2p - 1$ целых чисел содержит $p$ чисел, сумма которых делится на $p$.
\end{pr}

Этот забавный факт вытекает из следующего общего
результата (хотя может быть доказан и независимо).

\begin{pr}[Теорема Коши-Давенпорта]\label{t:cada} Пусть $p$ простое и $A,B\subset\Z_p$ непусты.
Обозначим
$A+B = \{a + b \,:\, a \in A,\; b \in B\}$.
Тогда $|A + B| \geqslant \min (|A| + |B| - 1, p)$.
\end{pr}

\emph{Указание.} Обозначим
\[
f(x,y) = \prod_{c\in A+B}(x+y-c).
\]
Тогда $f(x,y) = 0$ при всех $x\in A$, $y\in B$.

\begin{pr}\label{t:cada1} Пусть $p$ простое и $A,B\subset\Z_p$ непусты.
Обозначим
$$A\dotplus B = \{a+b \,:\, a\in A,\; b\in B, \; a\ne b\}.$$

(a) (Гипотеза Эрдеша-Хайльбронна) Имеем $|A \dotplus A| \geqslant \min (2|A| - 3, p).$

(b) Если $A\ne B$, то $|A \dotplus B| \geqslant \min (|A| + |B| - 2, p)$.
\end{pr}

\emph{Указание.} П. (a) сводится к п. (b).
П. (b) сводится к случаю $|A|\ne|B|$, который доказывается при помощи
теоремы \ref{c:degb'}.b.

\begin{pr}\label{p:plan} (a) (Международная олимпиада 2007)
Найдите наименьшее целое $k$, такое что в пространстве $\R^3$ имеется $k$ плоскостей,
объединение которых не содержит точку $(0,0,0)$, но содержит все остальные точки $(a,b,c)\in\Z^3$, такие что $0\le a,b,c\le n$.


(b) (теорема Олсона) Пусть $p$ простое.
Найдите наименьшее целое $m$, для которого среди любых $m$ векторов из $\Z_p^k$ (не обязательно различных) найдутся несколько векторов  с нулевой суммой.
\end{pr}

\begin{pr}\label{p:reg} Пусть $p$ простое, степень каждой вершины графа меньше $2p$, а средняя степень
больше $2p-2$.
Тогда имеется непустой
подграф, степень каждой вершины которого равна $p$.
\end{pr}


\begin{pr}\label{p:snevily} Пусть $k$ и $n$~--- целые положительные числа.

(a) Если $2k\le n+1$, то для любых $a_1,\ldots,a_k \in \Z_n$ (не обязательно различных) существует перестановка $\sigma$ множества $\{1, \ldots, k\}$ такая, что вычеты $a_1 + \sigma(1), \ldots, a_k + \sigma(k) \in \Z_n$ попарно различны.

(b) Найдите коэффициент при одночлене $x_1^{k - 1}\cdots x_k^{k - 1}$  многочлена
\[
\prod_{1\le i<j\le n} (x_j-x_i)^2.
\]

(c) (гипотеза Сневиля, доказанная Алоном)
Пусть $p$~--- нечётное простое число, $k<p$ и вычеты $b_1,\ldots,b_k\in\Z_p$ попарно различны.
Тогда для любых $a_1, \ldots, a_k \in\Z_p$ (не обязательно различных) найдётся перестановка $\sigma$ множества  $\{1,\ldots,k\}$ такая, что вычеты $a_1+b_{\sigma(1)},\ldots,a_k+b_{\sigma(k)} \in \Z_p$ попарно различны.
\end{pr}

Здесь п. (b) --- подсказка к п. (a), а п. (c) --- обобщение п. (a) для нечётного простого $n$.

\begin{pr}\label{p:diff} Каждое из данных $2^n+1$ конечных множеств покрашено в один из двух цветов.
Оба цвета имеются.
Тогда найдётся не менее $2^n$ попарно различных множеств (не обязательно из данных), каждое из которых является симметрической разностью
двух множеств разных цветов.
(Напомним, что \textit{симметрическая разность} двух множеств есть множество всех элементов, принадлежащих ровно одному из них.)
\end{pr}

\bigskip
{\bf Решения задач}

\smallskip
{\it В решениях задач \ref{p:setp}, \ref{t:cada}, \ref{t:cada1},
\ref{p:plan}.b, \ref{p:reg}
равенства являются равенствами вычетов по модулю $p$ или многочленов с коэффициентами в $\Z_p$.}

\smallskip
{\bf \ref{p:setp}.}
Если среди данных $2p - 1$ чисел найдутся $p$ чисел, дающих одинаковые остатки от деления на $p$, то возьмем их.
Иначе упорядочим остатки $x_1 \geq x_2 \geq \ldots \geq x_{2p - 1}$ от деления данных чисел на $p$.
Обозначим $M_p = \{x_p\}$ и $M_i = \{x_i, x_{p + i}\}$ для $i = 1,\ldots, p - 1$.
Индукцией по $n$ с помощью теоремы Коши-Давенпорта \ref{t:cada} нетрудно получить, что
\[
|A_1 + \ldots + A_n| \geq \min(|A_1| + \ldots + |A_n| - n + 1, p)
\]
для любых непустых подмножеств $A_1,\ldots,A_n \subset \Z_p$.
Поэтому $|M_1 + \ldots + M_p| = p$.
Значит, $0 \in M_1 + \ldots + M_p$.

\smallskip
{\bf \ref{t:cada}.} \textit{Набросок доказательства,
не использующего теоремы \ref{c:degb'}.}
Можно считать, что $|A|\leq |B|$.
Проведём индукцию по $|A|$.
База $|A| = 1$ очевидна.

Переход. Пусть $|A|>1$. Если $B = \Z_p$, то теорема очевидна.
Иначе, ввиду простоты числа $p$, существует такое $x\in\Z_p$, что для $A' := A+x$ выполнено $A'\cap B\neq \varnothing$ и $A'\not\subset B$.
Тогда по индукции можно перейти к паре множеств $(A'\cap B, A'\cup B)$.

\smallskip
{\bf \ref{t:cada}.}
Обозначим $a:=|A|$, $b:=|B|$, $C := A + B$.

Пусть сначала $a+b-1>p$.
Тогда существуют непустые подмножества $A'\subset A$ и $B'\subset B$, для которых $|A'|+|B'| = p+1$.
Поскольку $A'+B' \subset C$ и $\min(a+b-1,p) = p = \min(|A'|+|B'|-1, p)$, то утверждение сведено к случаю  $a+b-1=p$.

Пусть теперь\footnote{Аналогично предыдущему абзацу утверждение сводится к случаю $a+b-1=p$.
Вместо этого мы включаем в нижеприведенную формулу множитель $(x+y)^{a+b-2-|C|}$.}
$a+b-1\le p$.
Пусть, напротив, $|C|\le a+b-2$.
Положим
\[
f(x,y) = (x+y)^{a+b-2-|C|} \prod_{c\in C}(x+y-c).
\]
Ясно, что  $f(x,y) = 0$ при всех $x\in A$, $y\in B$.
Тогда по теореме \ref{c:degb'}.b коэффициент при $x^{a-1}y^{b-1}$ равен нулю.
Но он равен ${a+b-2\choose a-1}$.
Так как $a+b-2\le p-1$, то он  не делится на $p$.
Противоречие.

\smallskip
{\bf \ref{t:cada1}.}
(a)  Возьмём произвольный элемент $a\in A$.
Тогда $A \dotplus A = A\setminus\{a\} \dotplus A$, и мы свели задачу к пункту (b).

(b) Обозначим $a:=|A|$, $b:=|B|$, $C:=A\dotplus B$.
Если $a+b-2>p$, то аналогично доказательству теоремы Коши-Давенпорта \ref{t:cada} удалим из множеств $A$ или $B$ какие-нибудь элементы и сведем утверждение к случаю  $a+b-2=p$.

Пусть теперь $a+b-2\le p$.
Пусть, напротив, $|C|\le a+b-3$.
Положим
\[
f(x,y) = (x-y) (x+y)^{a+b-3-|C|} \prod_{c\in C}(x+y-c).
\]
Ясно, что $f(x,y) = 0$ при всех $x\in A$, $y\in B$.
Тогда по теореме \ref{c:degb'}.b коэффициент при $x^{a-1}y^{b-1}$ равен нулю.
Но он равен ${a+b-3\choose a - 2} - {a+b-3\choose a - 1}$.

Так как $a+b-3 \le p-1$, то при $a\neq b$ эта разность не делится на $p$.
Противоречие.
(Этого достаточно для п. (а).)

Пусть теперь $a=b$.
Если $A\cap B = \varnothing$, то применим теорему Коши-Давенпорта \ref{t:cada}.
Иначе перейдём ко множествам $(A', B') = (A\cap B, A\cup B)$.
Нетрудно видеть, что
$$A' + B' \subset A + B,\quad |A'| + |B'| - 2 = |A| + |B| - 2\quad\text{и}\quad |A'| < |B'|.$$
Таким образом, утверждение сведено к случаю $a\ne b$.

\smallskip
{\bf \ref{p:plan}.} (a) Ответ: $3n$.

В качестве примера возьмем
\begin{itemize}
\item $n$ плоскостей $x=a$, $1\leq a\leq n$;

\item $n$ плоскостей $y=b$, $1\leq b\leq n$;

\item $n$ плоскостей $z=c$, $1\leq c\leq n$.
\end{itemize}

Пусть теперь $k<3n$ и множество $\{0,1,\ldots,n\}^3-\{(0,0,0)\}$ покрыто $k$ плоскостями с уравнениями $a_ix+b_iy+c_iz+ d_i = 0$.
Обозначим
\[
f(x,y,z) = \prod_{i=1}^k(a_ix+b_iy+c_iz+ d_i).
\]
Тогда $\deg f=k<3n$.
Этот многочлен равен нулю во всех точках множества $\{0,1,\ldots,n\}^3$, кроме, быть может, точки $(0,0,0)$.
Тогда по теореме \ref{c:degb}.b он равен нулю и в точке $(0,0,0)$.
Значит, некоторая плоскость содержит точку $(0,0,0)$.

(b)  Ответ: $k(p - 1) + 1$.

Ясно, что $m > k(p - 1)$, поскольку можно выбрать базисные векторы $e_1$, $e_2$, $\ldots$, $e_{k}$ по $p - 1$ раз.

Рассмотрим произвольные $m = k(p - 1) + 1$ векторов $v_1, \ldots, v_m \in \mathbb{Z}_p^k$.
Покажем, что можно выбрать несколько из них с нулевой суммой.
Сумма нескольких из этих векторов~--- линейная комбинация $\varepsilon_1 v_1 +  \ldots + \varepsilon_m v_m$, где $\varepsilon_i \in \{0, 1\}$ для $i \in \{1, \ldots, m\}$.
Обозначим $j$-ю координату вектора $v_i$ через $v_{ij}$  для $i \in \{1, \ldots, m\}$, $j \in \{1, \ldots, k\}$.
Если $\varepsilon_1 v_1 +  \ldots + \varepsilon_m v_m \neq 0$, то существует такое $j \in \{1, \ldots, m\}$, что $\varepsilon_1 v_{1j} + \ldots + \varepsilon_m v_{mj} \neq 0$.
В пространстве $\Z_p^k$ с координатами $\varepsilon_1,\ldots,\varepsilon_m$ рассмотрим $m-1=k(p-1)$ гиперплоскостей, заданных уравнениями
\[
\varepsilon_1v_{1j} + \ldots + \varepsilon_m v_{mj} = \lambda, \quad \text{где}\quad
i \in \{1, \ldots, k\}\quad \text{и}\quad \lambda \in \Z_p\setminus\{0\}.
\]
Аналогично версии п. (a) по модулю $p$ эти гиперплоскости не могут покрывать все точки куба $\{0,1\}^m$, кроме точки $(0,\ldots,0)$.
Значит, существует точка $(\varepsilon_1,\ldots,\varepsilon_m) \in \{0, 1\}^m$, не лежащая ни в одной из построенных $k(p-1)$ гиперплоскостей.
Тогда $\varepsilon_1v_1 + \ldots + \varepsilon_m v_m = 0$.

\smallskip
{\bf \ref{p:reg}.} Обозначим через $V$ и $E$ множества вершин и рёбер данного графа.
Для каждого ребра $e\in E$ введем переменную $x_e$.
Определим многочлен с коэффициентами в $\Z_p$ формулой
\[
f(x) = \prod_{v\in V}\left(\left(\sum_{e\in E\,:\,v\in e}x_e\right)^{p-1}-1\right).
\]
Имеем $f(0,\ldots,0)\neq 0$.
Кроме того, $\deg f\leq|V|(p-1)<|E|$
(если граф имеет изолированную вершину, то в определении многочлена $f$ суммирование ведется по пустому множеству, поэтому степень многочлена $f$ не равна $|V|(p - 1)$).
Значит, по теореме \ref{c:degb}.b многочлен принимает ненулевое значение ещё хотя бы в одной точке $z\in \{0,1\}^{|E|}$.
Удалим из графа все те рёбра $e$, для которых $z_e=1$.
Если появились изолированные вершины, то удалим их.
Получим искомый подграф.


\smallskip
{\bf \ref{p:snevily}.}
(a) Для $x_i, x_j \in \{1, \ldots, k\}$ имеем $|x_i - x_j| \leq k - 1 \leq (n - 1)/2 < n/2$ и
\[
x_i + a_i \not\equiv x_j + a_j\,\,(\operatorname{mod}\,n) \quad\Leftrightarrow\quad x_j - x_i \not\equiv a_i - a_j\,\,(\operatorname{mod}\,n) \quad\Leftrightarrow\quad x_j - x_i \neq r_{ij},
\]
где $r_{ij}$ обозначает число, сравнимое с $a_i - a_j$ по модулю $n$, из интервала $(-n/2, n/2]$.
Таким образом, достаточно показать, что существуют различные $x_1,\ldots, x_k \in \{1, \ldots, k\}$ такие, что $x_j - x_i \neq r_{ij}$ для любых $1 \leq i < j \leq k$.
Ввиду теоремы \ref{c:degb'}.b достаточно показать, что коэффициент при мономе $x_1^{k-1}\cdots x_k^{k-1}$ многочлена
\[
\prod_{1 \leq i < j \leq k} (x_j - x_i)(x_j - x_i - r_{ij})
\]
не равен нулю.
Так как этот коэффициент совпадает с коэффициентом при том же мономе многочлена
\[
\prod_{1 \leq i < j \leq k} (x_j - x_i)^2,
\]
то нужное вытекает из ответа на п. (b).

(b) Ответ: $k! (-1)^{{k\choose2}}$.

Для доказательства воспользуемся формулой (для определителя Вандермонда)
\[
\prod_{1\le i<j\le k} (x_j-x_i) =
\sum_{\sigma} \operatorname{sgn}\sigma \cdot x_{\sigma(1)}^{0}\cdots x_{\sigma(k)}^{k-1},
\]
где $\sigma$ пробегает все перестановки множества $\{1, \ldots, k\}$, а $\operatorname{sgn}\sigma$~--- знак перестановки.
Тогда
\[
\prod_{1\le i<j\le k} (x_j - x_i)^2 =
\left(\sum_{\sigma} \operatorname{sgn} \sigma \cdot x_{\sigma(1)}^{0}\cdots x_{\sigma(k)}^{k - 1}\right)
\cdot \left(\sum_{\pi} \operatorname{sgn} \pi \cdot x_{\pi(1)}^{0}\cdots x_{\pi(k)}^{k - 1}\right).
\]
При раскрытии скобок моном $x_1^{k - 1}\cdots x_k^{k - 1}$ присутствует только в произведении слагаемых, соответствующих всем тем перестановкам $\sigma$ и $\pi$, для которых $\pi = \sigma\alpha$, где $\alpha(j)=k+1-j$ при любом $j$.
Для таких перестановок
$\operatorname{sgn} \sigma \cdot \operatorname{sgn} \pi = \operatorname{sgn} \alpha = (-1)^{{k\choose2}}$.
Поэтому искомый коэффициент равен $k!(-1)^{{k\choose2}}$.

\bigskip
{\bf Литература}

[A] N. Alon, Combinatorial Nullstellensatz, Combin. Probab. Comput. 8 (1999), 7--29,

\url{https://www.cs.tau.ac.il/~nogaa/PDFS/null2.pdf}

[BS] A. Bishnoi and P.L. Clark, Restricted variable Chevalley-Warning theorems,

\url{http://alpha.math.uga.edu/~pete/Bishnoi-Clark22.pdf}

\smallskip
[C] E. Chen, Combinatorial Nullstellensatz, 2013

\url{https://web.evanchen.cc/handouts/BMC_Combo_Null/BMC_Combo_Null.pdf}

\smallskip
[C14] P.L. Clark, The Combinatorial Nullstellens\"atze Revisited, The Electr. J. of Comb., 21 (2014),
\url{https://www.combinatorics.org/ojs/index.php/eljc/article/view/v21i4p15}

\smallskip
[CGS] P.L. Clark, T. Genao, and F. Saia, Chevalley-Warning at the Boundary, Expositiones Math., to appear,

\url{http://alpha.math.uga.edu/~pete/Chevalley_Warning_on_the_Boundary.pdf}

\smallskip
[D] В. Димитров. Combinatorial Nullstellensatz, стр. 416-425 в книге
<<Петербургские олимпиады школьников по математике 2003-2005>>.
СПб: Невский Диалект, БХВ-Петербург, 2007.

\smallskip
[KP] Р. Н. Карасёв, Ф. В. Петров. Ещё раз о комбинаторной теореме о нулях, стр. 116-120 в книге
<<Задачи Санкт-Петербургской олимпиады по математике 2010 года>>.
СПб: Невский Диалект, БХВ-Петербург, 2007.

\smallskip
[KS] A. Kezdy, H. Snevily, Distinct Sums Modulo $n$ and Tree Embeddings,  Combin. Probab. Comput. 11 (2002), 35--42,

\url{https://www.math.ucdavis.edu/~deloera/MISC/LA-BIBLIO/trunk/Kezdy/kezdy.pdf}

\smallskip
[P] Ф. В. Петров. Восстановление многочлена по его значениям, стр. 98-106 в книге
<<Задачи Санкт-Петербургской олимпиады школьников по математике 2015 года>>.
М.: МЦНМО, 2015.

\smallskip
[R1] {\it А.М. Райгородский,} Проблема Борсука, М.: МЦНМО, 2015.

\smallskip
[R2] {\it А.М. Райгородский,} Линейно-алгебраический метод в комбинаторике, М.: МЦНМО, 2015.

\smallskip
[S] {\it  А. Скопенков,} Короткое опровержение гипотезы Борсука, Мат. Просвещение, 17 (2013), 88--92,
\url{http://arxiv.org/abs/0712.4009}

\end{document}